\documentclass[10pt,a4paper]{article}
\usepackage{amsfonts,amssymb,amsmath,graphicx}
\usepackage{theorem} 
\usepackage{bm}      
\usepackage{eucal}   
\usepackage[english]{babel}

  {\theorembodyfont{\rm}
  \newtheorem{defi}{Definition}[section]

   }

  {\theorembodyfont{\it} 

    }
{\theorembodyfont{\rm}\theoremstyle{change}\newtheorem{nrtxt}[defi]{\hspace*{-0.55ex}}
\newenvironment{txt}{\begin{trivlist}\item}{\end{trivlist}}

\newlength\rest 

\let\phi=\varphi

\let\theta=\vartheta

\newcommand{\bC}{{\mathbb C}}
\newcommand{\bD}{{\mathbb D}}

\newcommand{\bP}{{\mathbb P}}

\newcommand{\bR}{{\mathbb R}}

\newcommand{\bZ}{{\mathbb Z}}

\newcommand{\cC}{{\mathcal C}}

\newcommand{\cG}{{\mathcal G}}

\newcommand{\cL}{{\mathcal L}}
\newcommand{\cM}{{\mathcal M}}
\newcommand{\cN}{{\mathcal N}}

\DeclareMathOperator{\diag}{diag}

\DeclareMathOperator{\GE}{GE} 
\DeclareMathOperator{\GL}{GL}

 \DeclareMathOperator{\Trans}{T}

\newcommand{\dis}{\mathbin{\scriptstyle\triangle}}
\newcommand{\notdis}{\mathbin{\not\scriptstyle\triangle}}
\newcommand{\Matrixfeld}[4]{\left#1\!\begin{array}{*{#3}{c}}#4\end{array}\!\right#2}

\newcommand{\Mat}{\Matrixfeld()}

\begin{document}
\title{From pentacyclic coordinates to chain geometries, and back\thanks{This article
is an extended version of my lecture ``\emph{Jenseits der
penta\-zyklischen Koordinaten}'' at the colloquium celebrating the 75th
birthday of Prof.~Dr.~Dr.~h.~c. Walter Benz, University of Hamburg, June 9th,
2006.}}
\author{Hans Havlicek\\Institut f\"{u}r Diskrete Mathematik und Geometrie,\\
        Technische Universit\"{a}t Wien,\\
        Wiedner Hauptstra{\ss}e 8--10,\\
        A-1040 Wien, Austria,\\
        \mbox{\tt havlicek@geometrie.tuwien.ac.at}}
\date{}

\maketitle

\begin{abstract}\noindent
Starting with the classical circle geometry of Sophus Lie, we give a survey
about some of the developments in the area of chain geometries during the last
three decades.
\par\noindent
MSC 2000: 51B05, 51B10, 51B15, 51B20, 51B25, 51A45, 51C05, 17C50, 14M15.
\par\noindent
Key words: Lie geometry, projective line over a ring, distant graph,
Grassmannian, chain geometry, Jordan isomorphism, Jordan system.
\end{abstract}

\section{Pentacyclic Coordinates}

\begin{nrtxt}
The circle geometry of S.~Lie (1842--1899) aims at eliminating the distinction
between circles, lines, and points of the Euclidean plane. The idea is that
points and lines should be viewed as circles with ``zero'' and ``infinite''
radius, respectively. Moreover, lines and circles are endowed with an
\emph{orientation}. For our purposes it suffices to think of an \emph{oriented
line\/} (an \emph{oriented circle\/}) as a line (a circle) with an arrow on it.
There are precisely two possibilities of orientation. For circles we can even
distinguish between \emph{counterclockwise\/} and \emph{clockwise\/}
orientation. More precisely, a \emph{Lie cycle\/} is one of the following:
\par
\begin{itemize}
\item
An oriented circle. Its \emph{signed radius\/} $r\neq 0$ is positive (negative)
if the orientation is counterclockwise (clockwise).
\item A point. Its \emph{radius} is defined to be zero.
\item An oriented line.
\item The \emph{point at infinity\/}. It is denoted by the symbol
$\infty$.
\end{itemize}
The set of all Lie cycles will be written as $\cN$. It is endowed with a binary
\emph{contact relation\/} $\sim$, where $x\sim y$ is to be read as ``$x$
touches $y$''. This relation is reflexive and symmetric by definition. Touching
Lie cycles are depicted in Figure~\ref{abb:kontakt}. In addition, $\infty$ is
assumed to touch all oriented lines, but no oriented circle and no point.\\

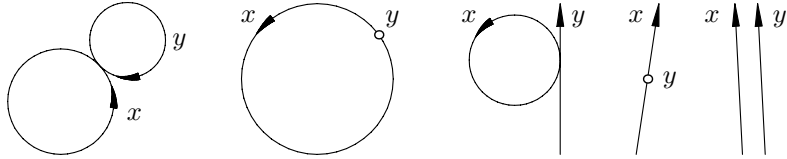
\begin{figure}[h!]\unitlength1cm
\centering
\begin{picture}(10,2.2)
    \put(0,0){\includegraphics[width=20cm]{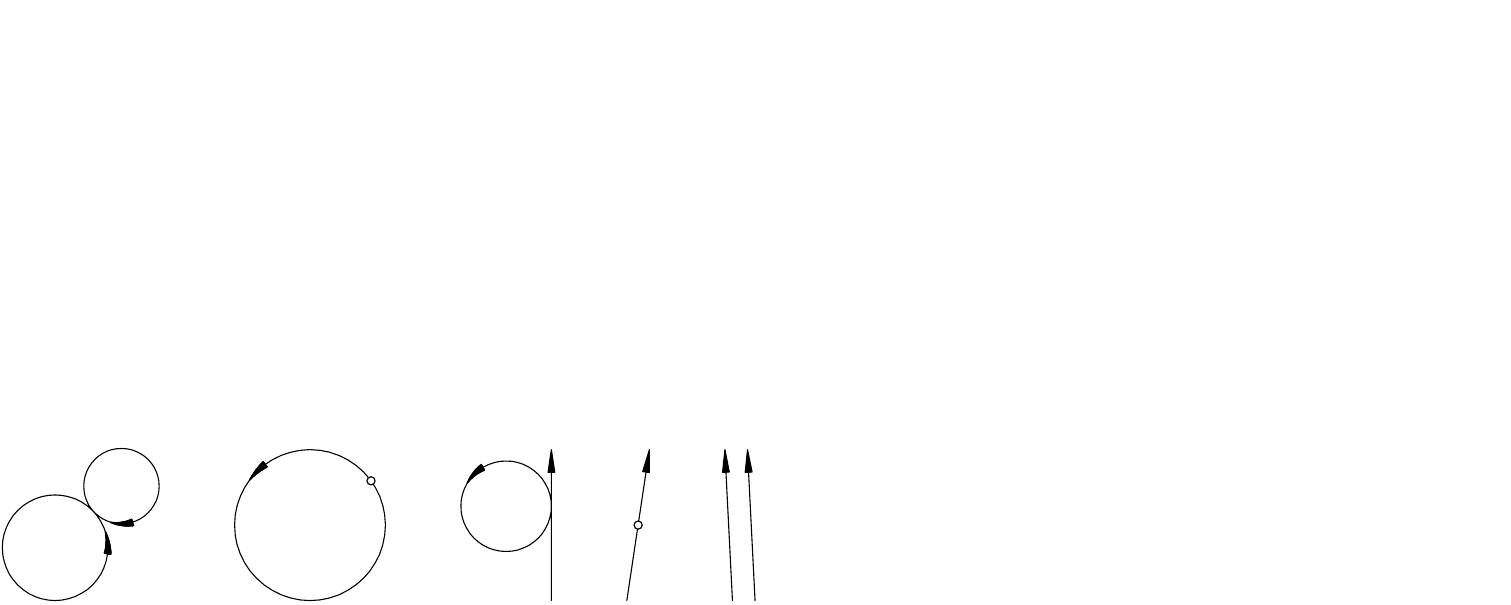}}
    \put(1.6,0.5){$x$}\put(2.2,1.5){$y$}
    \put(3.1,1.8){$x$}\put(5.0,1.8){$y$}
    \put(6.0,1.8){$x$}\put(7.45,1.8){$y$}
    \put(8.2,1.8){$x$}\put(8.65,1.0){$y$}
    \put(9.2,1.8){$x$}\put(10.1,1.8){$y$}
\end{picture}
\caption{Contact relation}\label{abb:kontakt}
\end{figure}\vspace{-\belowdisplayskip}
\end{nrtxt}

\begin{nrtxt}
Let us shortly motivate the ``need'' for orientation and for the point at
infinity: Suppose that we are given three distinct points. There is a
\emph{unique\/} circle or a \emph{unique\/} line passing through all of them.
On the other hand, if we are given the three side lines of a triangle then
there are \emph{four\/} circles touching all of them, the incircle and the
three excircles. Thus points and (non-oriented) lines behave totally different.
How should they be considered as being ``equal''? However, our three given
points give rise to precisely \emph{two\/} distinct Lie cycles touching all of
them. So, let us also introduce an orientation on each of the three given
lines. Then precisely \emph{one\/} of the four circles from the above can be
oriented in such a way that it is in contact with the oriented lines. This is
better than before, but we would like to have \emph{two\/} such circles. Hence
we add an extra \emph{point at infinity\/} which touches all spears
irrespective of their orientation.
\end{nrtxt}

\begin{nrtxt}
The set $\cN$ of Lie cycles can be mapped bijectively onto the point set of a
non-degenerate quadric
\begin{equation}\label{eq:Lie}
    \Lambda: -x_0^2+x_1^2+x_2^2+x_3^2-x_4^2=0
\end{equation}
in the four-dimensional real projective space $\bP_4(\bR)$ as follows: Choose a
Cartesian coordinate system in the plane.
\begin{itemize}

\item
The image of an oriented circle with midpoint $(m_1,m_2)$ and signed radius
$r\neq 0$ is
    $$\bR\left(\frac{1+N}2,\frac{1-N}2,m_1,m_2,-r\right),$$
where $N:=m_1^2+m_2^2-r^2$. The image of a point $(m_1,m_2)$ is given likewise
by setting $r=0$.

\item
In order to obtain the image point of an oriented line we consider its two
equations in Hesse normal form. Precisely one of them, say
$a_0+a_1x_1+a_2x_2=0$, has the property that the unit vector $(-a_2,a_1)$
determines the orientation of the given line. The image point is then defined
to be
    $\bR\left(-a_0,a_0,a_1,a_2,1 \right).$

\item The image of $\infty$ is $\bR(-1,1,0,0,0)$.
\end{itemize}
The \emph{pentacyclic coordinates\/} of a Lie cycle are, by definition, the
homogeneous coordinates of its image point on the \emph{Lie quadric\/}
$\Lambda$.
\par
It is easy to check that two Lie cycles are in contact if, and only if, their
images on $\Lambda$ are conjugate with respect to the polarity of $\Lambda$ or,
in other words, if their pentacyclic coordinate vectors are orthogonal with
respect to the pseudo-Euclidean dot product of $\bR^5$ given by the matrix
$\diag(-1,1,1,1,-1)$.
\par
A \emph{Lie transformation\/} is a bijection $\cN\to\cN$ which preserves $\sim$
in both directions. The \emph{fundamental theorem\/} of Lie geometry states
that all Lie transformations arise from the collineations of $\bP_4(\bR)$
leaving invariant the Lie quadric $\Lambda$ \cite[p.~42]{benz-73}.
\end{nrtxt}

\begin{nrtxt}
In order to illustrate the power of the mapping described in the above, we
recall a problem due to Apollonius of Perga (approx.\ 262--190 B.~C.):
\emph{Find all circles which touch three given circles.}
\par
In terms of pentacyclic coordinates the solution can be found most easily:
First, endow each of the given circles with an orientation. Next, calculate
their images on the Lie quadric $\Lambda$. Then intersect the tangent
hyperplanes of $\Lambda$ at these points. This gives (up to degenerate cases) a
line $L$, say. Finally, determine $L\cap\Lambda$. This amounts to solving a
quadratic equation, whereas all steps before yield linear equations. If the
line $L$ and the Lie quadric $\Lambda$ have points in common (which need not be
the case) then their pre-images are the solutions to the Apollonius problem for
oriented circles. See Figure~\ref{abb:apollonius}. Taking into account the
various possibilities for orientation of the given circles, the initial problem
turns out to have up to eight solutions. By this approach, also all forms of
exceptional cases, e.~g., when one of the ``solutions'' is an oriented line or
a point, are easily understood.
\begin{figure}[h!]\unitlength1cm
\centering
\begin{picture}(3.2,3.4)
    \includegraphics[width=20cm]{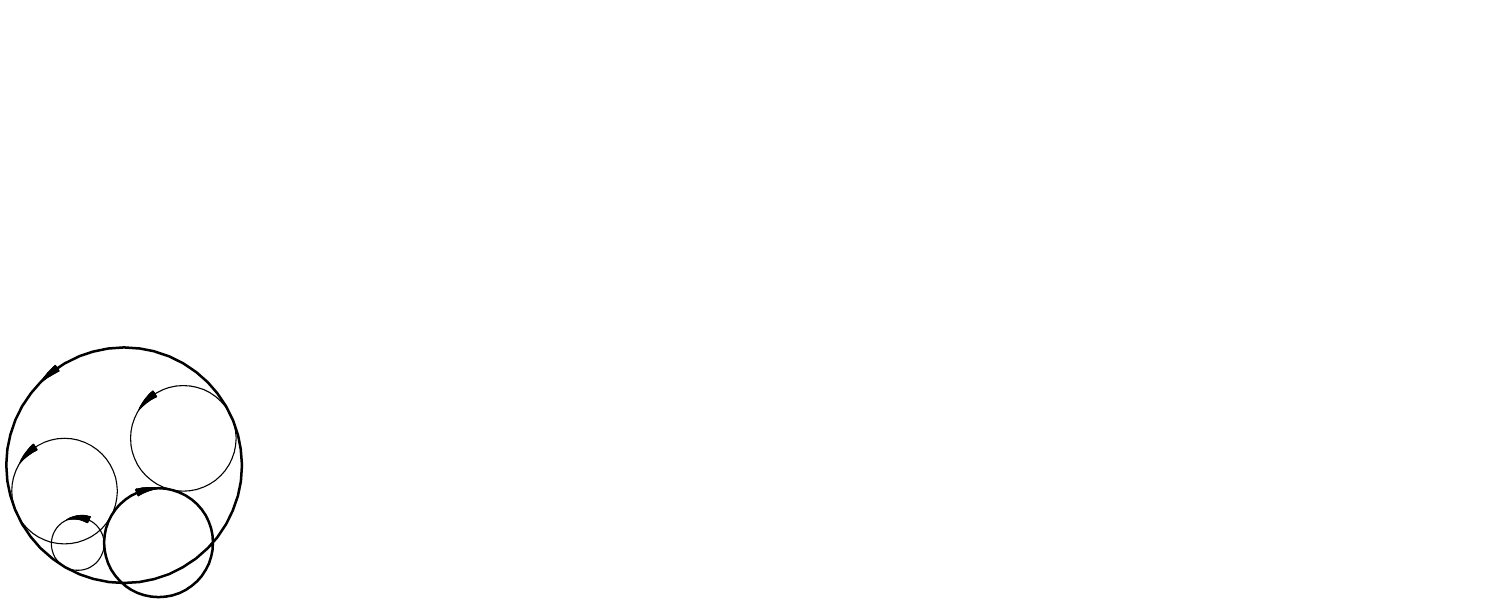}
\end{picture}
\caption{The problem of Apollonius for oriented circles}\label{abb:apollonius}
\end{figure}\vspace{-\belowdisplayskip}
\end{nrtxt}

\begin{nrtxt}
We refer to \cite{benz-73} for more details on Lie's circle geometry and an
extensive list of references up to the year 1973. It is worth noting that the
point-line geometry of the Lie quadric is one of the classical
\emph{generalised quadrangles} \cite[pp.~57--58]{vanmaldeg-98}. For
higher-dimensional Lie geometries, differential Lie geometry, and relations to
special relativity see \cite{benz-92}, \cite{cecil-92}, \cite{evel+mc+t-74},
\cite{fill-79}, \cite{giering-82}, \cite{karz+k-88}, \cite{metzger-87},
\cite{rigby-81}, \cite{rigby-81a}, \cite{yaglom-81}, and the references made
there. Infinite-dimensional Lie geometry is one of the topics in the recent
book of W.~Benz \cite{benz-05}.
\end{nrtxt}

\begin{nrtxt}
We now aim at recovering two other classical geometries from Lie's circle
geometry.
\par
Firstly, let $\cM$ be the set of points in the plane together with $\infty$.
For each Lie cycle $y\notin\cM$ the point set
\begin{equation}\label{eq:moebiuskreise}
    \{x\in\cM\mid x\sim y\}
\end{equation}
is the set of points on a circle or the range of points (including $\infty$) on
a line. This set remains unchanged if the orientation of $y$ is altered. In
this way we obtain the Euclidean \emph{M\"{o}bius geometry\/}
({August~F.~M\"{o}bius} (1790--1868)). In contrast to Lie geometry, which is a
set endowed with a binary relation, here we have a set of points together with
the family of distinguished subsets (\ref{eq:moebiuskreise}) carrying the name
\emph{M\"{o}bius circles\/} or \emph{chains}. Observe that M\"{o}bius circles
do not have an orientation.
\par
Secondly, let $\cL$ be the set of oriented lines which now will also be called
\emph{spears}. Each Lie cycle $y\notin\cL\cup\{\infty\}$ gives rise to the set
\begin{equation}\label{eq:speerkette}
    \{x\in\cL\mid x\sim y\}
\end{equation}
which is called a \emph{chain of spears\/}, shortly a \emph{chain}. This gives
the Euclidean \emph{Laguerre geometry}, i.~e.\ the set $\cL$ together with the
set of all its chains. It is named after Edmond~N.~Laguerre (1834--1886). The
point $\infty$ is superfluous in Laguerre geometry.
\par
Both geometries allow a unified description. If we consider the usual field of
\emph{complex numbers}
then the point set of the Euclidean M\"{o}bius geometry coincides with the
\emph{complex projective line},
\begin{equation*}
    \bP(\bC):=\bC\cup\{\infty\}\mbox{~with~}\infty:=\frac10,
\end{equation*}
which is well known from complex analysis. Likewise, the ring of \emph{real
dual numbers},
\begin{equation}\label{eq:dualzahlen}
    \bD:=\{x+ y\varepsilon \mid
    (x,y)\in\bR^2\}, \mbox{~where~}\varepsilon\notin\bR\mbox{~and~}\varepsilon^2=0,
\end{equation}
gives rise to the \emph{dual projective line}. It has the form
\begin{equation}\label{eq:P-dual}
    \bP(\bD):=\bD\cup\left\{\frac{1}{y\varepsilon}\mid y\in\bR\right\}.
\end{equation}
In contrast to the complex projective line there are \emph{infinitely many}
points at infinity; they are given by the second set on the right hand side of
(\ref{eq:P-dual}). This set comprises all formal quotients with enumerator $1$
and a zero divisor as denominator. By a classical result, the set of spears of
the Euclidean Laguerre geometry can be identified with the dual projective line
\cite[pp.~26--28]{benz-73}.
\par
In either case the chains are precisely the images of the \emph{real projective
line\/} $\bP(\bR):=\bR\cup\{\infty\}$, considered as subset of $\bP(\bC)$
(resp.\ $\bP(\bD)$), under the action of the complex (resp.\ dual) linear
fractional group
\begin{equation*}
    z\mapsto \frac{az+b}{cz+d} \mbox{~~with~}\Mat2{a & b\\c &d}\mbox{~invertible.}
\end{equation*}
\end{nrtxt}

\begin{nrtxt}
In the famous book \cite{benz-73} by W.~Benz, published in 1973, projective
lines over \emph{commutative\/} rings and chain geometries arising from
\emph{commutative\/} algebras were investigated systematically, thereby
generalising the classical results. In this article we focus our attention on
the further development of these topics.
\par
The book of Benz contains a wealth of further material which we cannot mention
here.
\end{nrtxt}

\section{The projective line over a ring}

\begin{nrtxt}
We adopt the following conventions: All our rings are associative, with a unit
element $1\neq 0$, which acts unitally on modules, and is inherited by
subrings. Multiplication in a \emph{field\/} need not be commutative. If a ring
$R$ contains a field $K$, as a subring which commutes with all elements of $R$,
then $R$ is called a \emph{$K$-algebra}. The dimension of $R$ over $K$ may be
finite or infinite.
\end{nrtxt}

\begin{nrtxt}
The crucial task is to find a ``good'' definition of the \emph{projective line}
over a ring $R$, even when $R$ is not necessarily commutative. In terms of left
homogeneous coordinates a \emph{point} of this line should of course have the
form $R(a,b)$ with $(a,b)\in R^2$ (considered as a left module over $R$). But,
which pairs $(a,b)$ should be representatives of points? Let us shortly recall
some particular cases:
\begin{itemize}
\item
If $R$ is a field then $(a,b)$ gives rise to a point if, and only if $(a,b)\neq
(0,0)$. (This is mathematical folklore.)

\item
If $R$ is a \emph{local ring}, i.~e., the set of all non-invertible elements of
$R$ is an ideal, then $(a,b)$ determines a point if, and only if, $a$ or $b$ is
an invertible element.
(This was used, e.~g., in \cite{lima+l-77a}.)

\item
If $R$ is commutative then a pair $(a,b)$ yields a point if, and only if, it is
\emph{unimodular\/}, i.~e., there are elements $x,y\in R$ with $ax+by=1$. (This
definition was adopted in \cite{benz-73}.)
\end{itemize}
All these conditions remain meaningful over any ring. Hence there are
\emph{different definitions\/} for the projective line over a ring. See
\cite{lash-97} for a survey and \cite[pp.~291--292]{veld-85} for further
comments on the problem of obtaining ``good'' projective geometries from a
ring.
\end{nrtxt}

\begin{nrtxt}
Projective lines over several classes of non-commutative rings, like \emph{skew
fields} (see \cite{benz-73}), \emph{matrix rings over commutative fields}, and
rings of \emph{ternions} (i.~e.\ upper triangular $2\times 2$ matrices over a
commutative field), were exhibited already in the 1960s and before. The
Belgians J.~Depunt, \cite{depu-59}, \cite{depu-60}, C.~Vanhelleputte
\cite{vanhelle-66}, {X.\ Hubaut} \cite{hubaut-64}, \cite{hubaut-65}, and
{J.~A.~Thas} \cite{thas-69}, \cite{thas-71} were among the first to study
projective lines over non-commutative rings other than skew fields. The Italian
G.~Russo considered the projective line over a ring of upper triangular
$m\times m$ matrices over a commutative field. He coined the name ``ennoni''
(Italian for ``$n$-ions'') for such a ring \cite{russo-65}, \cite{russo-67},
\cite{russo-72}.
\par
We shall stick here to a definition which, to our knowledge, appeared first in
\cite{hubaut-65}. There is even a footnote in this article pointing out the
general case of an arbitrary ring $R$, whereas the paper itself is concerned
with finite-dimensional algebras only. The essential ingredient for this
definition is the \emph{general linear group\/} $\GL_2(R)$ in two variables
over a ring $R$. The elements of this group are precisely the invertible
$2\times 2$ matrices with entries in $R$.
\end{nrtxt}

\begin{defi}
The \emph{projective line\/} over a ring $R$ is the set $\bP(R)$ of all cyclic
submodules $R(a,b)$ of $R^2$, where $(a,b)$ is the first row of an invertible
$2\times 2$ matrix over $R$. Such a pair is called \emph{admissible}.
\end{defi}

\begin{txt}
We read off from the $2\times 2$ identity matrix that $R(1,0)$ is a point.
Furthermore, as $A$ ranges in $\GL_2(R)$ all points of $\bP(R)$ arise as
$R(1,0)\cdot A$. Thus $\bP(R)$ can also be described as the \emph{orbit\/} of
$R(1,0)$ under the natural action of $\GL_2(R)$. This means that all points of
$\bP(R)$ are ``the same'' up to the action of $\GL_2(R)$. It is an easy
exercise to show that two admissible pairs represent the same point if, and
only if, they are left-proportional by an invertible element in $R$.
\par
An elegant coordinate-free definition of $\bP(R)$ is due to A.~Herzer; see
\cite[p.~785]{herz-95}. We refer also to the recently published book by
A.~Herzer and A.~Blunck \cite{blunck+he-05} for further details.
\end{txt}

\begin{nrtxt}\label{merkwuerdig}
The projective line over a ring $R$ has some peculiar properties if there exist
elements with a single-sided inverse \cite{blunck+h-00b}:
\par
If $s\in R$ has a left-inverse, say $l$, such that $ls=1\neq sl$ then
$R(s,0)=R(ls,0)=R(1,0)$, but it is easily seen that there is no matrix in
$\GL_2(R)$ with first row $(s,0)$. Thus a point may have non-admissible
representatives. However, from now on \emph{only admissible pairs will be used
to represent points}.
\par
If $s\in R$ has a right inverse, say $r$, such that $sr=1\neq rs$ then
$R(s,0)\subsetneqq R(1,0)$, but now $(s,0)$ is admissible, because
\begin{equation*}
    \gamma:=\Mat2{s&0\\1-rs&r}
    \mbox{~~has the inverse~~}  \gamma^{-1}=\Mat2{r&1-rs\\0&s}.
\end{equation*}
This means that there may be \emph{nested points\/}. So, we are far away from
Euclid's definition: ``\emph{A point is that which has no part}''
\cite[Vol.~1,~p.~153]{euclid-56}.
\par
There is another phenomenon which is only present in certain non-commutative
rings. \emph{A unimodular pair need not be admissible}
\cite[Remark~5.1]{blunck+h-01a}. The following example is based on a paper of
M.~Ojanguren and R.~Sridharan \cite{ojan+s-71}: Let $K$ be a skew field and let
$R=K[X,Y]$ be the polynomial ring in two independent central indeterminates
over $K$. Given elements $a,b\in K$ with $ab\neq ba$ the pair $(X+a,Y+b)$ turns
out to be unimodular, but not admissible.
\par
On the other hand, each admissible pair $(a,b)$ is unimodular, as follows by
multiplying an invertible matrix $A={a\;\;b\choose *\;\;*}$ with its inverse. A
ring has \emph{stable rank $2$} if for each unimodular pair $(a,b)\in R^2$
there is a $c\in R$ such that $a+bc$ is invertible. For such a ring unimodular
and admissible pairs are the same \cite[p.~785]{herz-95}. F.~D.~Veldkamp
(1931--1999) has repeatedly stressed the importance of rings of stable rank $2$
for geometry; see, e.~g., \cite[p.~293]{veld-85}.

\end{nrtxt} 

\section{The distant graph}

\begin{nrtxt}
On the projective line over a ring $R$ there is an important binary relation:
Two points $p$ and $q$ of $\bP(R)$ are called \emph{distant\/}, in symbols
$p\dis q$, if there is a matrix
\begin{equation*}
    \Mat2{a & b\\c & d}\in\GL_2(R)
\end{equation*}
with $p=R(a,b)$ and $q=R(c,d)$. This relation is anti-reflexive and symmetric.
Distant points are also said to form a \emph{clear\/}, \emph{spectral\/} or
\emph{regular\/} pair \cite{hotje-76}, \cite{hubaut-65}. The graph of the
relation $\dis$, i.~e. the pair $\big(\bP(R),{\dis}\big)$, is called the
\emph{distant graph\/} of $\bP(R)$. It is an undirected graph without loops.
\par
Other authors prefer the negated relation $\notdis$ and speak of
\emph{neighbouring\/} or \emph{parallel\/} points. The term ``parallel'' stems
from parallel spears in Euclidean Laguerre geometry, as depicted on the right
hand side of Figure~\ref{abb:kontakt}. Parallel spears correspond to
non-distant points of the dual projective line. The parallelism of spears is an
equivalence relation, but $\notdis$ is \emph{not\/} transitive in general.
\end{nrtxt}

\begin{nrtxt}\label{beispiele}
We take a closer look at some examples.
\begin{enumerate}
\item\label{beispiele.a}
Let $\bZ_4$ be the \emph{ring of integers modulo\/} $4$. We consider also the
ring of \emph{dual numbers over the field $\bZ_2$}. These dual numbers are
defined as in (\ref{eq:dualzahlen}), with $\bR$ to be replaced by $\bZ_2$. The
distant graphs of $\bP(\bZ_4)$ and $\bP(\bZ_2[\varepsilon])$ are isomorphic to
the graph of vertices and edges of an octahedron
(Figure~\ref{abb:distgraphen}). But our two rings are not isomorphic, since we
have $1+1\neq 0$ in $\bZ_4$ and $1+1=0$ in $\bZ_2[\varepsilon]$. Thus
non-isomorphic rings may have isomorphic distant graphs.

\item\label{beispiele.b}
The ring $\bZ_2\times\bZ_2$ of \emph{double numbers\/} over $\bZ_2$ has also
four elements. The distant graph of $\bP(\bZ_2\times\bZ_2)$ is depicted in
Figure~\ref{abb:distgraphen}, the shaded triangles serve only for better
visualisation. There are nine points.
\begin{figure}[h!]
\begin{minipage}\leftmargin
${}$ 
\end{minipage}
\setlength\rest\textwidth%
\addtolength\rest{-\leftmargin}%
\begin{minipage}\rest
\centering\unitlength3mm
\begin{picture}(15.67,7.76)
    \put(0,1.05){\includegraphics[width=5.65\unitlength]{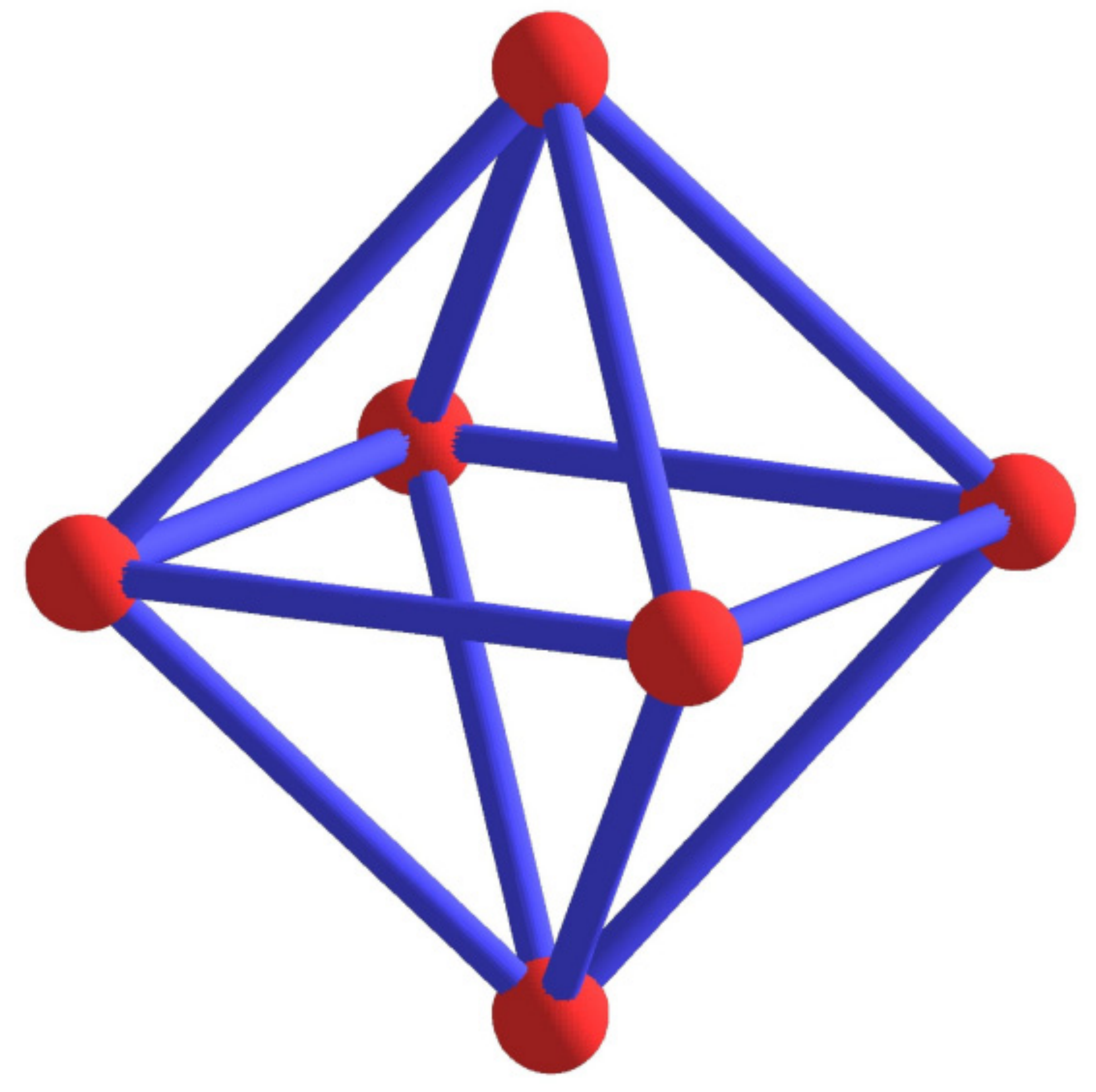}}
    \put(10,0){\includegraphics[width=5.67\unitlength]{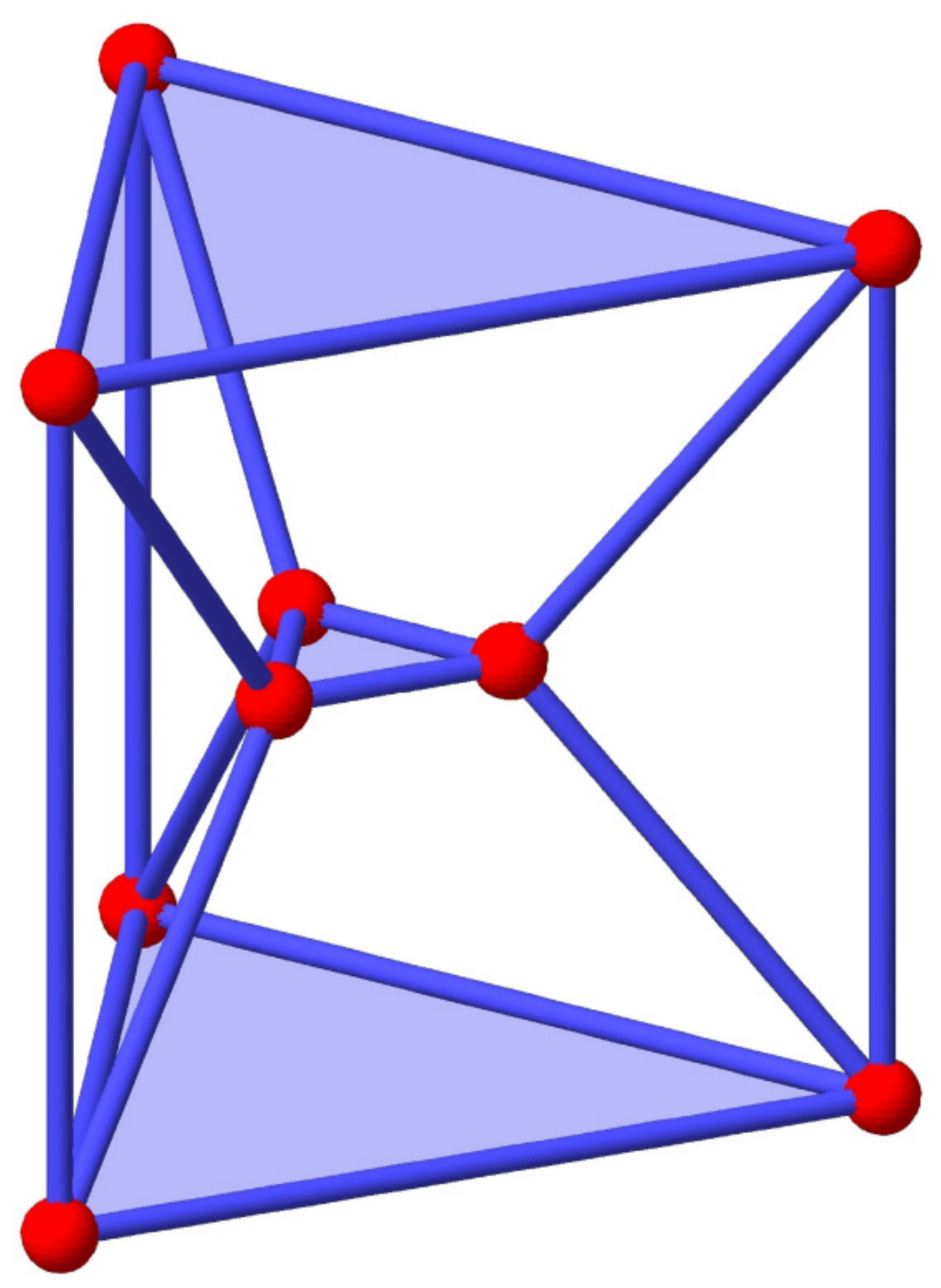}}
\end{picture}
    \caption{Distant graphs for rings with four elements.~}\label{abb:distgraphen}
\end{minipage}
\end{figure}\vspace{-\belowdisplayskip}
\item
The projective line over the field $\bC$ of complex numbers can be seen as the
unit sphere in Euclidean space. Here and in the subsequent examples it will be
more intuitive to illustrate the neighbour relation $\notdis$. Each point
$p\in\bP(\bC)$ has a single neighbour, namely $p$ itself
(Figure~\ref{abb:nachbarn}).
\begin{figure}[ht!]\unitlength3mm
\begin{minipage}\leftmargin
${}$ 
\end{minipage}
\setlength\rest\textwidth%
\addtolength\rest{-\leftmargin}%
\begin{minipage}\rest
\centering%
\begin{picture}(29.76,6.76)
    \put(0,0.76){\includegraphics[width=5.40\unitlength]{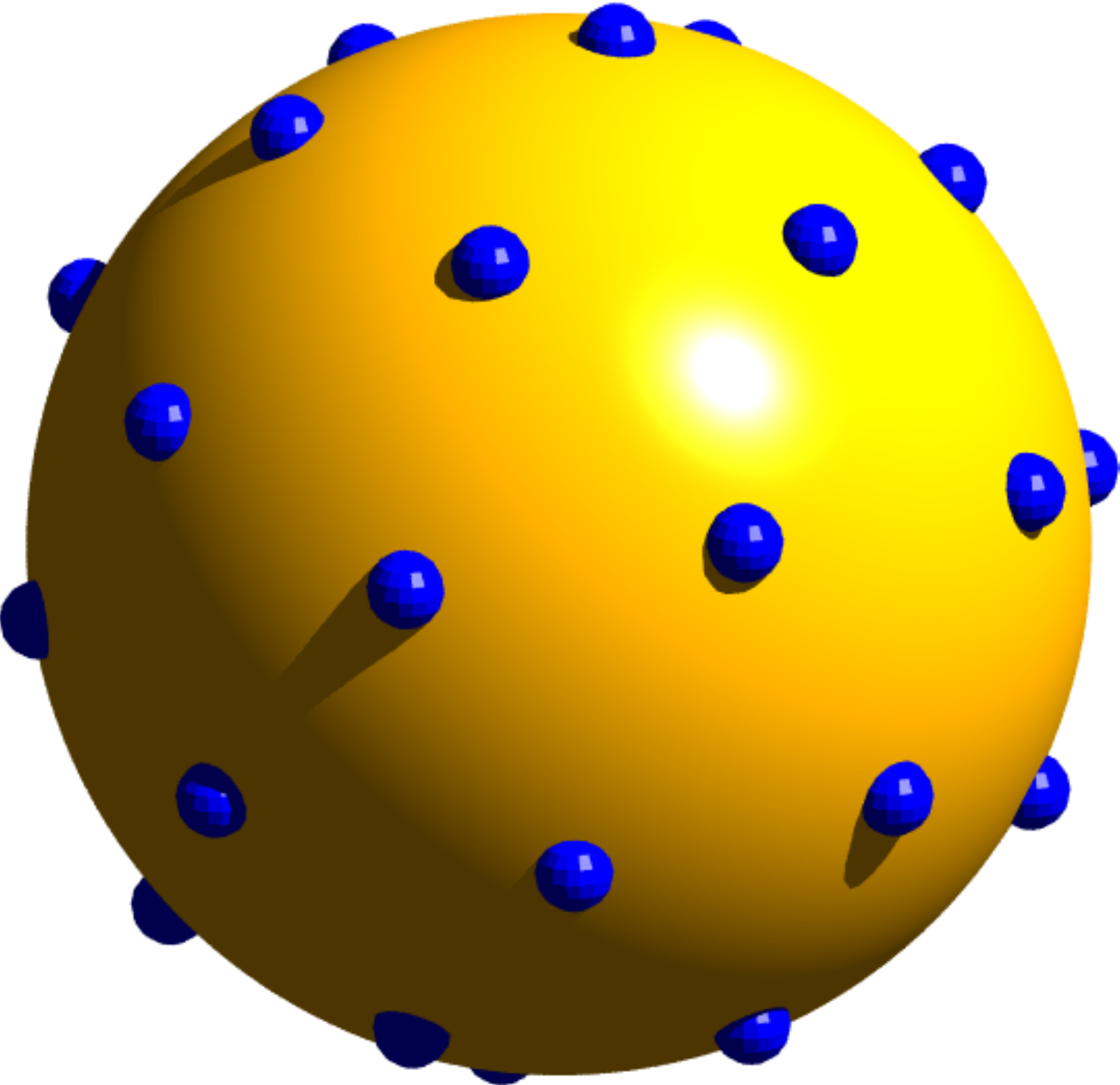}}
    \put(10.0,1.37){\includegraphics[width=8.40\unitlength]{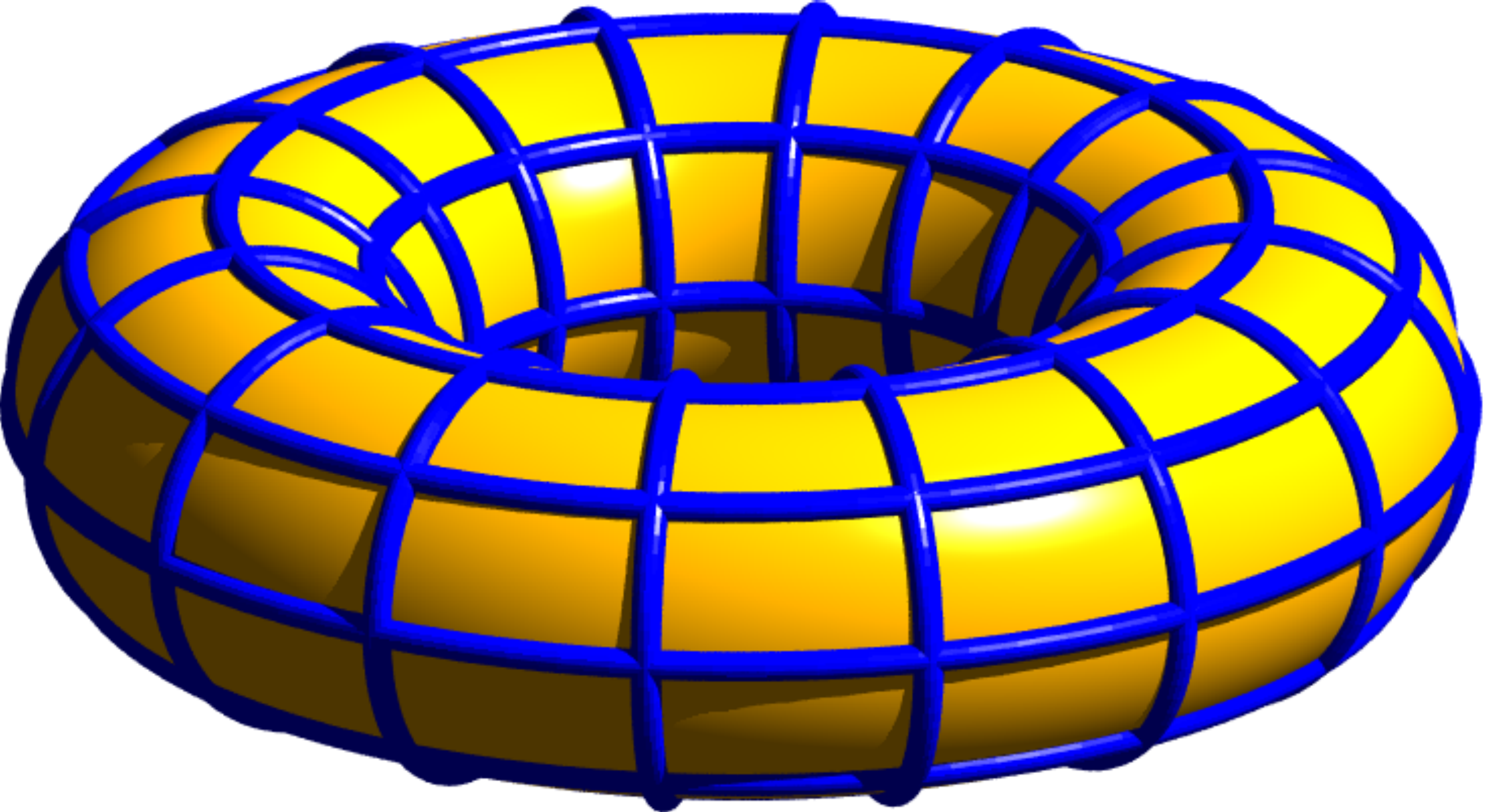}}
    \put(23,0){\includegraphics[width=5.26\unitlength]{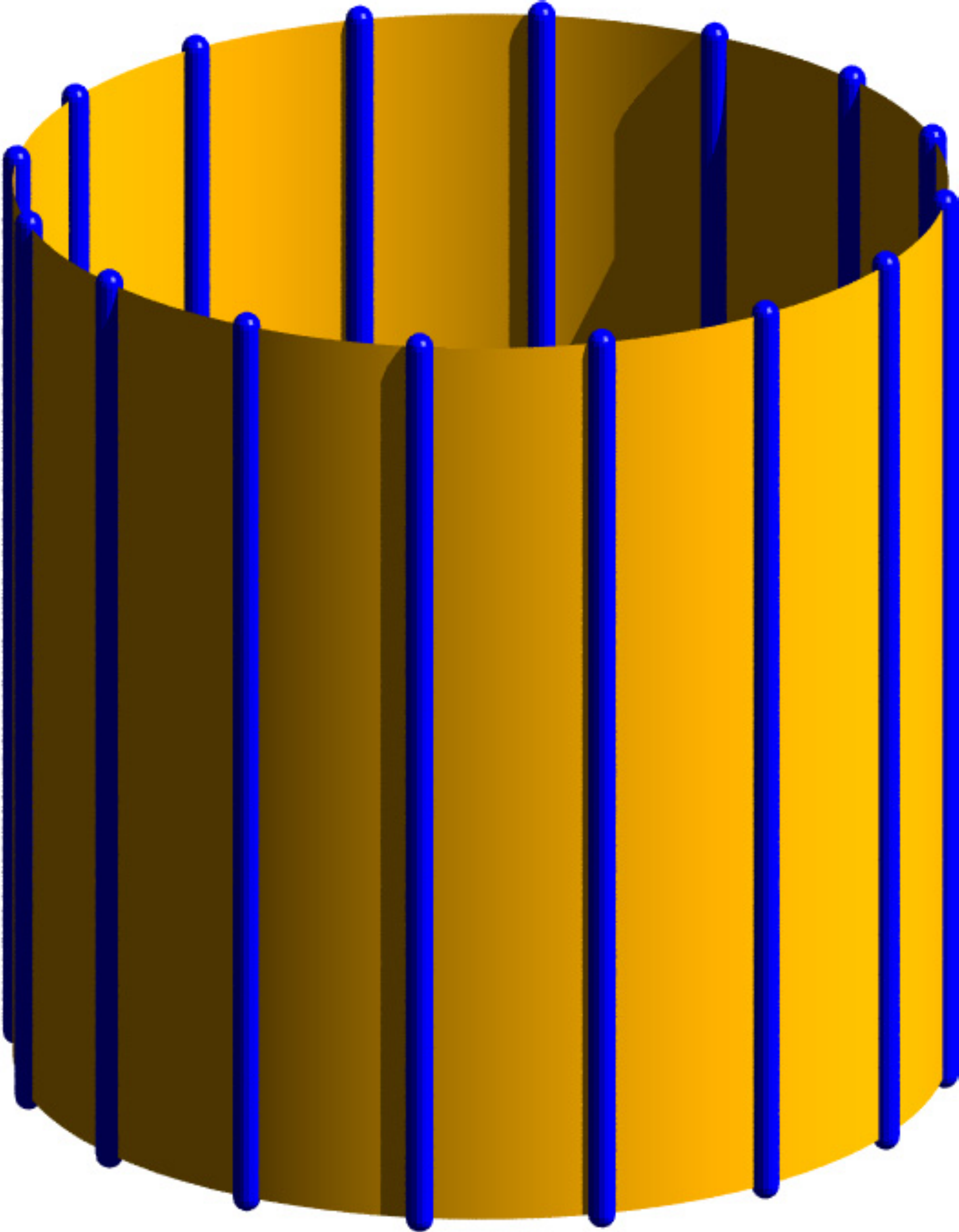}}
\end{picture}%
\caption{Neighbour relation on $\bP(\bC)$, $\bP(\bR\times\bR)$, and
    $\bP(\bD)$.~}\label{abb:nachbarn}%
\end{minipage}
\end{figure}\vspace{-\belowdisplayskip}

\item
The projective line over the ring $\bR\times \bR$ of \emph{real double
numbers\/} can be identified with the Cartesian product
$\bP(\bR)\times\bP(\bR)$. By virtue of this identification, two points
$p=(p_1,p_2)$ and $q=(q_1,q_2)$ with $p_i,q_j\in\bP(\bR)$ are neighbouring if
$p_1=q_1$ or $p_2=q_2$. Since the real projective line may be illustrated as a
circle, the torus is a point model for the projective line over $\bR\times\bR$.
The longitudinal and latitudinal circles represent the maximal subsets of
mutually neighbouring points (Figure~\ref{abb:nachbarn}).

\item
The \emph{Blaschke cylinder}, named after W.~Blaschke (1885--1962), is a point
model for the projective line over the real dual numbers. The generators of the
cylinder represent the maximal subsets of mutually neighbouring points
(Figure~\ref{abb:nachbarn}).

\item\label{bsp:matrizenring}
Let $K^{n\times n}$ be the ring of $n\times n$ matrices over a commutative
field $K$. The projective line over $K^{n\times n}$ can be identified with the
Grassmannian $\cG_n$ of all $n$-dimensional subspaces of $K^{2n}$ via the
bijection
\begin{equation}\label{eq:grassmann}
    \bP(K^{n\times n})\to \cG_n : R(A,B)\mapsto \mbox{rowspace of~}(A|B),
\end{equation}
where $(A|B)$ stands for the matrix $A$ augmented by $B$. Under this bijection
distant points correspond to complementary subspaces. See
\cite[p.~500]{hubaut-65}, where the result is stated in an equivalent form,
using the projective space $\bP_{2n-1}(K)$. We know today from the work of
A.~Blunck that this theorem holds also for matrix rings over skew fields and,
mutatis mutandis, for \emph{endomorphism rings of infinite-dimensional vector
spaces} over arbitrary fields \cite[Theorem~2.1]{blunck-99}.
\end{enumerate}
\end{nrtxt}

\begin{nrtxt}
Among all examples from above the one in (\ref{bsp:matrizenring}) is most
important. By the results of E.~Artin (1898--1962) and J.~H.~M.~Wedderburn
(1882--1948), the Artinian semisimple rings are precisely the direct products
of matrix rings over (possibly different) fields \cite{lam-91}. The projective
line over such a ring is in one-one correspondence with a \emph{product of
Grassmannians\/} \cite{blunck+h-05b}. This gives a deep insight into the
structure of such a projective line. For example, this allows one to determine
the number of points of the projective line over a finite ring, even if the
ring is not semisimple \cite[pp.~31--36]{veld-89}.
\end{nrtxt}

\begin{nrtxt}
Let us turn back to the general case. The distant graph of a projective line
over a ring $R$ is easily seen to be a complete graph if, and only if, $R$ is a
field. In this case the diameter of the distant graph equals one. By
\cite{herz-95} this diameter is $\leq 2$ whenever $R$ is a ring of stable rank
$2$. For example, finite-dimensional algebras are of stable rank $2$. This
explains why in most ``classical examples'' for any two non-distant points
there is a point which is distant to both of them. Only few distant graphs with
a diameter $>2$ seem to be known \cite{blunck+h-01a}.
\par
The distant graph of $\bP(R)$ is connected if, and only if, $R$ is a
\emph{$\GE_2$-ring\/}. This means that each matrix in $\GL_2(R)$ is a product
of elementary matrices and invertible diagonal matrices. The results of
{P.~M.~Cohn} \cite{cohn-66} provide examples of distant graphs with more than
one connected component \cite[p.~115]{blunck+h-01a}.
\par
Only recently, the meaning of the \emph{Jacobson radical\/} of a ring was
expressed in terms of the associated distant graph \cite[p.~116]{blunck+h-03a}:
If $p=R(1,0)$ and $q\in\bP(R)$ then
\begin{equation*}
    (x\dis p \Leftrightarrow x\dis q) \mbox{~for all~}x\in\bP(R)
\end{equation*}
holds precisely when $q=R(1,r)$ with $r$ taken from the Jacobson radical of
$R$.
\end{nrtxt}

\begin{nrtxt}
The group $\GL_2(R)$ acts transitively on the set of mutually distant triples
of $\bP(R)$ and thereby preserves the distant relation. It is a natural
question to ask for all isomorphisms between distant graphs. In
\cite{blunck+h-05a} an equivalent problem for Grassmannians of vector spaces is
exhibited, whereas in \cite{blunck+h-05b} all isomorphisms $\bP(R)\to\bP(R')$
are determined provided that $R$ and $R'$ are direct products of matrix rings
over fields. In both papers there are many interrelations with results about
\emph{adjacency preserving transformations} of Grassmannians due to {W.~Benz},
{H.~Brauner}, {W.-l.~Huang}, {A.~Kreuzer}, {A.~Naumowicz}, {K.~Pra\.{z}mowski},
and others. See \cite{benz-92}, \cite{brau-88}, \cite{havl-95},
\cite{huang-98}, \cite{kreuz-98}, and \cite{nau+p-01}. The distant graph
$\big(\bP(R),{\dis}\big)$ turns into a \emph{Pl\"{u}cker space\/}
\cite[p.~199]{benz-92} if we add a loop at each point. The \emph{Pl\"{u}cker
transformations} are then just the automorphisms of the distant graph.
\end{nrtxt}

\section{Chain geometries}

\begin{nrtxt}
Throughout this section $R$ denotes an algebra over a (necessarily commutative)
field $K$. The definition of the \emph{chain geometry} associated with $K$ and
$R$ can be taken over literally from \cite{benz-73}, since we did already
introduce the projective line over $R$ irrespective of commutativity: As $R$
has a unit element, we may assume that $K\subset R$. This allows to identify
$\bP(K)$ with a subset of $\bP(R)$ via $K(a,b)\mapsto R(a,b)$. Every image of
$\bP(K)$ under a matrix of $\GL_2(R)$ is called a \emph{chain\/}. Let $\cC$ be
the set of all chains. Then the incidence structure
\begin{equation*}
    \Sigma(K,R):=\big(\bP(R),\cC\big)
\end{equation*}
is called the \emph{chain geometry\/} associated with the $K$-algebra $R$.
\par
The distant relation on $\bP(R)$ can be expressed in terms of $\Sigma(K,R)$ as
follows: Two distinct points are distant if, and only if, they are on a common
chain. Thus all results about the distant graph are also available in a chain
geometry.
\end{nrtxt}

\begin{nrtxt}
Many basic properties of a chain geometry do not depend on the commutativity of
the algebra $R$. For example, there is a unique chain through any three
mutually distant points. A major difference between commutative and
non-commutative algebras (or, more generally, rings) concerns the notion of
\emph{cross ratio\/} \cite[p.~787]{herz-95}. In the non-commutative case the
cross ratio of four points of $\bP(R)$ is a class of conjugate elements of $R$
rather than a single element of $R$. Cross ratios are often useful, since four
mutually distinct points are on a common chain precisely when their cross ratio
is in $K$. We refer in particular to the paper of A.~Blunck \cite{blunck-03a}
dealing with cross ratios on Grassmannians.
\end{nrtxt}

\begin{nrtxt}\label{punktmodelle}
It is well known that the M\"{o}bius geometry $\Sigma(\bR,\bC)$, the Laguerre
geometry $\Sigma(\bR,\bD)$, and the \emph{Minkowski geometry\/}
$\Sigma(\bR,\bR\times\bR)$ (named after H.~Minkowski (1864--1909)) can be
represented on an elliptic quadric, a quadratic cone (without its vertex), and
a hyperbolic quadric in three-dimensional real projective space $\bP_3(\bR)$.
In either case the non-degenerate conics represent the chains. For the
M\"{o}bius geometry and Laguerre geometry one may also use a Euclidean space
and view the elliptic quadric and the cone as a sphere and a cylinder of
revolution (Figure~\ref{abb:nachbarn}). However, the hyperbolic quadric
modelling the Minkowski geometry cannot be seen as part of the
three-dimensional Euclidean space. This is why we used a torus in
Figure~\ref{abb:nachbarn} instead. The grid of longitudinal and latitudinal
circles corresponds to the grid of lines on the hyperbolic quadric.
\end{nrtxt}

\begin{nrtxt}
One of the open problems from \cite{benz-73} was to find \emph{point models\/}
for arbitrary chain geometries. In view of the examples from
\ref{punktmodelle}, it was quite natural to look for \emph{quadric models},
where chains were represented by non-degenerate conics. H.~Hotje showed in
\cite{hotje-74} and \cite{hotje-76} that such models exist for
finite-dimensional \emph{quadratic algebras}, i.~e.\ $K$-algebras in which
every element has a minimal polynomial over $K$ with degree $\leq 2$. Here the
work on these algebras by H.~Karzel (who used the term \emph{kinematic
algebras\/} \cite{karz-74}, \cite{karz+k-85}) turned out useful. Several papers
on this topic appeared afterwards. The next step was taken in 1980 by W.~Benz,
H.-J.~Samaga, and H.~Schaeffer \cite{benz+s+s-81}: The chain geometry
$\Sigma(K,K^n)$ was shown to be embeddable in a projective space over $K$ as
\emph{Segre variety}, with chains going over to \emph{normal rational curves}.
The breakthrough was accomplished shortly afterwards by M.~Werner
\cite{werner-82} using ideas from \cite{hubaut-65}. We sketch here Werner's
approach in a generalised form due to A.~Herzer \cite{herz-84a}:
\begin{itemize}
\item
Given a finite-dimensional $K$-algebra $R$, determine a \emph{faithful
representation\/} of $R$ in terms of $n\times n$ matrices over $K$. This
amounts to embedding $\Sigma(K,R)$ in $\Sigma(K,K^{n\times n})$ such that
chains of $\Sigma(K,R)$ go over to chains of $\Sigma(K,K^{n\times n})$.
\item
The bijection (\ref{eq:grassmann}) of $\bP(K^{n\times n})$ onto the
Grassmannian of $n$-subspaces of $K^{2n}$ gives a model of $\Sigma(K,R)$ within
this Grassmannian. Making use of results by {R.~Metz} \cite{metz-81}, the
images of chains can be identified as \emph{reguli\/} of the Grassmannian.
\item
Finally, this Grassmannian is mapped onto its associated Grassmann variety,
lying in a $\left({2n\choose n}-1\right)$-dimensional projective space. Here
the chains are represented by normal rational curves.
\end{itemize}
\par
All the classical point models mentioned in the above (excluding the torus) fit
into this general concept which may also be described in a coordinate-free way
\cite[p.~810]{herz-95}. By suitable projections, it is often possible to obtain
point models of $\Sigma(K,R)$ in lower-dimensional spaces. These ``projected
models'' are smooth \emph{quasiprojective varieties\/} and the chains appear
there as \emph{rational curves\/} \cite[p.~812]{herz-85a}.
\par
The case of infinite-dimensional algebras seems to be unsettled.
\end{nrtxt}

\section{Isomorphisms of chain geometries}

\begin{nrtxt}
Let $\Sigma(K,R)$ and $\Sigma(K',R')$ be chain geometries. It is one of the
basic problems to determine \emph{all\/} isomorphisms between them. In
\cite{benz-73} this problem was solved for various classes of chain geometries
over commutative algebras. However, in order to find generalisations one first
has to find appropriate mappings of the underlying algebras which can be used
to describe all isomorphisms in a second step.
\end{nrtxt}

\begin{nrtxt}\label{anti}
We start with an obvious example: Let $\alpha: R\to R'$ be an
\emph{isomorphism} of the $K$-algebra $R$ onto the $K'$-algebra $R'$ or, said
differently, let $\alpha$ be a semilinear bijection (of vector spaces) such
that $(ab)^\alpha=a^\alpha b^\alpha$ for all $a,b\in R$. It is obvious that
\begin{equation}\label{1}
    R(a,b) \mapsto R'(a^\alpha,b^\alpha)
\end{equation}
defines an \emph{isomorphism\/} of $\Sigma(K,R)$ onto $\Sigma(K',R')$.
\par
Since we admit non-commutative algebras, we may also consider an
\emph{antiisomorphism\/} of algebras, i.~e.\ a semilinear bijection $\alpha:
R\to R'$ such that $(ab)^\alpha=b^\alpha a^\alpha$ for all $a,b\in R$. It was a
longstanding open problem whether or not any antiisomorphism defines ``in some
natural way'' an isomorphism of chain geometries. Observe that the assignment
given by (\ref{1}) is not well-defined in this case, let alone its being an
isomorphism: For if $u$ is a unit in $R$ then $(a,b)$ and $(ua,ub)$ represent
the same point, whereas $(a^\alpha,b^\alpha)$ and $(a^\alpha u^\alpha,b^\alpha
u^\alpha)$ need not be left proportional by a unit in $R'$. Hence they may
represent distinct points.
\par
If we restrict ourselves to local algebras then each point of $\bP(R)$ has at
least one \emph{normalised\/} representative $(1,b)$ or $(a,1)$. Now its image
can be defined unambiguously by $R'(1',b^\alpha)$ or $R'(a^\alpha,1')$. This
gives not only a well defined mapping, but also an isomorphism of chain
geometries. It was shown in \cite{blunck+h-01b} that any antiisomorphism of
algebras gives rise to an isomorphism of the associated chain geometries.
However, there does not seem to be an explicit formula for this isomorphism in
the general case. We refrain from a further discussion, because there are even
more general mappings of algebras which deserve our attention!
\end{nrtxt}

\begin{nrtxt}
Isomorphisms and antiisomorphisms are just particular examples of \emph{Jordan
isomorphisms} of algebras (P.~Jordan (1902--1980)). A Jordan isomorphism
$\alpha: R\to R'$ is a semilinear bijection taking $1$ to $1'$ such that
\begin{equation*}
    (aba)^\alpha=a^\alpha b^\alpha a^\alpha \mbox{~for all~} a,b\in R.
\end{equation*}
\end{nrtxt}

\begin{nrtxt}
For many classes of algebras, e.~g.\ fields or commutative algebras over fields
of characteristic $\neq 2$, there are no Jordan isomorphisms other than
isomorphisms and antiisomorphisms. On the other hand, proper Jordan
isomorphisms (other than isomorphisms and antiisomorphisms) are easy to
construct: The mapping which sends each $A\in\bR^{2\times 2}$ to its transpose
$A^{\Trans}$ is an antiautomorphism of $\bR^{2\times 2}$, whence the mapping
$(A,B)\mapsto (A,B^{\Trans})$ is a Jordan automorphism of the direct product of
$\bR^{2\times 2}$ with itself. Further examples of Jordan isomorphisms can be
found in \cite[pp.~81--82]{blunck+he-05}.
\end{nrtxt}

\begin{nrtxt}
Jordan isomorphisms of finite-dimensional local algebras determine isomorphisms
of their chain geometries, as was shown by {A.~Herzer} \cite{herz-87a}, who
could make use of previous results by {B.~V.~Limaye} and {N.~B.~Limaye}
\cite{lima+l-77b}, \cite{lima+l-77a}, and \cite{lima+l-77c}. Like before, the
main problem is the definition of such a mapping. Under the given restrictions
this can be done as for antiisomorphisms in \ref{anti} by using normalised
representatives.
\par
In 1989, C.~Bartolone made a great step forward by introducing a completely new
idea \cite{bart-89}. If $R$ has stable rank $2$ (see \ref{merkwuerdig}) then
\begin{equation*}
    \bP(R) = \{ R(xy-1,x)  \mid x,y\in R\}.
\end{equation*}
This means that each point of $\bP(R)$ can be written (usually in various ways)
with the help of two \emph{parameters} $x,y\in R$. Now, somewhat surprisingly,
the assignment
\begin{equation*}
    R(xy-1,x)\mapsto R'(x^\alpha y^\alpha -1',x^\alpha) \mbox{~with~} x,y\in R
\end{equation*}
gives a well defined isomorphism of chain geometries for any Jordan isomorphism
$\alpha:R\to R'$. A generalisation to arbitrary algebras and an interpretation
of Bartolone's approach can be found in \cite{blunck+h-03}. However the
definition of the point to point mapping arising from a Jordan isomorphism is
too involved to be sketched here. We just want to emphasise that this mapping
is defined only on the connected component of $R(1,0)$ in the distant graph.
This connected component may, or may not, be the entire projective line
$\bP(R)$.
\par
\end{nrtxt}

\begin{nrtxt}
The mappings from the previous paragraphs together with the mappings induced by
$\GL_2(R')$ give now all isomorphisms of chain geometries provided that certain
assumptions on $R$ and $R'$ are made. The interested reader should consult
\cite[pp.~832--833]{herz-95}. Also, we would like to add that the results on
mappings determined by Jordan isomorphisms can be reformulated in a more
general form for \emph{Jordan homomorphisms\/} of rings.
\end{nrtxt}

\section{Subspaces of chain geometries}

\begin{nrtxt}
Another remarkable topic is the investigation of \emph{subspaces\/} of a chain
geometry $\Sigma(K,R)$. A \emph{subspace\/} has to be closed under chains, but
it has also to satisfy a number of extra conditions in order to exclude
degenerate cases. For a precise definition one needs a series of notions which
are not within the scope of this article. See \cite[pp.~59--60]{blunck+he-05}.
Hence we have to restrict ourselves to presenting some examples of subspaces
together with their interesting algebraic background.
\end{nrtxt}

\begin{nrtxt}
Given a subalgebra $S$ of a $K$-Algebra $R$ one would expect the projective
line over $S$ to be a subspace of $\Sigma(K,R)$. However, this will only be
true if we impose the following extra condition:
\begin{equation}\label{eq:invers}
    a\in S\mbox{~invertible in~} R \Rightarrow a^{-1}\in S.
\end{equation}
If condition (\ref{eq:invers}) is not satisfied for an element $a\in S$ then
the distant relation $\dis_S$ on $\bP(S)$ does not coincide with the
restriction to $\bP(S)$ of the distant relation $\dis_R$ coming from $\bP(R)$.
Indeed, we have $R(a,1)\dis_R R(0,1)$, but $S(a,1)\,{\notdis}_S\, S(0,1)$. For
example, let $R$ be the polynomial algebra $K[X]$, and let $S$ be its field of
fractions $K(X)$. Then $X$ is invertible in $K(X)$, but $X^{-1}\notin K[X]$.
\end{nrtxt}
\begin{nrtxt}
While any subalgebra satisfying the extra condition (\ref{eq:invers}) gives
rise to a subspace $\Sigma(K,R)$, there are also more general substructures of
$R$ with this property: A \emph{strong Jordan system\/} of $R$ is defined to be
a $K$-subspace $S$ of $R$ such that (i) $S$ contains the unit element $1\in R$,
(ii) $S$ satisfies condition (\ref{eq:invers}), and (iii) for each $x\in R$
more than half of the elements in the coset $x+K$ are invertible (this is
strongness).
\par
Each strong Jordan system $S$ is closed under the \emph{Jordan triple product},
i.~e., $aba\in S$ for all $a,b\in S$, but it need not be closed under
multiplication. If the characteristic of $K$ is $\neq 2$ then
$\frac12\,(ab+ba)\in S$ for all $a,b\in S$. See
\cite[pp.~61--63]{blunck+he-05}.
\par
In order to associate with $S$ a subset of $\bP(R)$ the idea from
\cite{bart-89} to describe points of $\bP(R)$ via parameters $x,y\in R$ is used
once more: To each strong Jordan system $S$ corresponds the point set
\begin{equation}\label{eq:JS}
    \bP(S):= \{R(xy-1,x)\mid x,y\in S\}
\end{equation}
which turns out to be a subspace of $\Sigma(K,R)$ \cite[p.~67]{blunck+he-05}.
\par
For a wide class of algebras all connected subspaces are of the form
(\ref{eq:JS}), up to a transformation in $\GL_2(R)$. These results are due to
H.-J.~Kroll \cite{kroll-91}, \cite{kroll-92b}, \cite{kroll-92a}, and A.~Herzer
\cite{herz-92}. Cf.\ also \cite[p.~69--72]{blunck+he-05}.

\end{nrtxt}
\begin{nrtxt}
We consider the following example. Let $\bR^{2\times 2}$ be the algebra of
$2\times 2$ matrices over $\bR$. Denote by $S$ the subset of all
\emph{symmetric $2\times 2$ matrices}. While $S$ is not closed under
multiplication, it is a strong Jordan system of $\bR^{2\times 2}$. The chain
geometry $\Sigma(\bR,\bR^{2\times 2})$ has a quadric model, namely the well
known \emph{Klein quadric\/} in $\bP_5(\bR)$ representing the lines of
$\bP_3(\bR)$ (or, equivalently, the $2$-subspaces of the vector space $\bR^4$)
\cite[p.~814]{herz-95}. The subspace $\bP(S)$ corresponds to a hyperplane
section of the Klein quadric. This hyperplane section can be identified via a
collineation with the Lie quadric $\Lambda$ from (\ref{eq:Lie}). This means
that Lie's circle geometry allows an alternative description as the subspace of
$\Sigma(\bR,\bR^{2\times 2})$ associated to $S$. The chains correspond to the
non-degenerate conics on the Lie quadric. Figure~\ref{abb:apollonius}
illustrates one such chain within the set $\cN$ of Lie cycles: All Lie cycles
that touch the two bold cycles form a chain. (Only three of these cycles are
actually drawn.) Furthermore, Lie cycles represent distant points if, and only
if, they do not touch.
\end{nrtxt}
\section{Further reading}

\begin{nrtxt}
There are a lot of topics which would deserve our attention. Among them are
\emph{chain geometries over Jordan systems\/} and the algebraic description of
their isomorphisms via \emph{isotopisms\/} \cite{blunck+he-05},
\emph{generalised chain geometries\/} $\Sigma(K,R)$, where $K$ need not be in
the centre of the ring $R$ \cite{blunck+h-00a}, and \emph{geometries of field
extensions\/} \cite{havl-93b}, \cite{havl-94b}, \cite{havl-94a},
\cite{maeu+m+n-80}, \cite{piec-94}. There is a widespread literature about
characterisations of chain geometries and related structures which we did not
even touch upon. The book \cite{blunck+he-05} and the survey article
\cite{herz-95} are indispensable sources on these and various other topics. For
results and references on \emph{topological circle planes\/} and other
\emph{projective geometries over rings\/} see \cite{bert-04},
\cite{polster+s-01}, \cite{stei-95}, and \cite{veld-95}. A series of papers
deals with \emph{fractals\/} in chain geometries \cite{artzy-92a},
\cite{artzy-92c}, \cite{artzy-92b}, \cite{sama-95c}. Certain \emph{finite chain
geometries\/} lead to \emph{designs\/} or \emph{divisible designs\/}
\cite{blunck+h+z-06z}, \cite{giese+h+s-05a}, \cite{havl-06a}, \cite{spera-92a},
\cite{spera-95}. Connections between chain geometries and the \emph{geometry of
matrices\/} can be found in \cite{havl+s-06}, \cite{huang-05z}, and
\cite{wan-96}. The projective lines over some small rings found attention in
quantum physics \cite{planat+s+k-06a}, \cite{saniga+p-06a},
\cite{saniga+p-06b}, \cite{saniga+p-06x}, \cite{saniga+p-06y},
\cite{saniga+p-06z}, \cite{saniga+p+k+p-06z}, \cite{saniga+p+m-06z},
\cite{saniga+p+p-06z}, \cite{saniga+p+p-06y}. For applications in twistor
theory of the projective point and the projective line over biquaternions see
\cite{agnew-03}, \cite{soucek-81}, and \cite{soucek-82}.
\end{nrtxt}

\begin{nrtxt}
Acknowledgement. The author is greatly indebted to Andrea Blunck (Hamburg),
Herbert Hotje (Hannover), and Armin Herzer (Bodman) for their support in
preparing this article.
\end{nrtxt}

\fontsize{9}{10pt}\selectfont

\end{document}